\documentstyle{article}
\begin{document}
\renewcommand{\thefootnote}{\fnsymbol{footnote}}

\begin{center}{\Large\bf BEYOND THE `PENTAGON IDENTITY'.}
\bigskip

{\large\bf A. Yu. Volkov\footnote{On leave of absence from
Steklov Mathematical Institute, St. Petersburg;\\
supported by INTAS grant CT93-0023}}

{\em Physique - Math\'{e}matique,
Universit\'{e} de Montpellier II\\
Pl. E. Bataillon, Case 50,
34095 Montpellier C\'{e}dex 05, France}
\end{center}
\bigskip
\begin{quote}{\bf Abstract.} An algebraical background of
the Lattice Conformal Field Theory is refined with the help
of a novel $q$-exponential identity.
\end{quote}
\bigskip

\noindent It is commonly believed [GR] that the function
$$ s(x)=\prod_{n\geq 0}(1-xq^{2n+1})                       $$
is a $q$-world counterpart of the exponential function.
It means that as soon as $u$ and $v$ make a Weyl pair
$$ uv=q^2 vu                                               $$
the $q$-exponents of them behave just like ordinary exponents
of commuting arguments do:
$$ s(u)s(v)=s(u+v).                                        $$
Recently [FV] added a missing `reversed' multiplication rule
$$ s(v)s(u)=s(u+v-qvu)                                     $$
to the collection of its properties. This time let me present
another identity
$$ s(v)s(u^{-1})s(u)s(v)=s(u^{-1})s(v)s(u)                 $$
which is a consequence of the two multiplication
rules but apparently has virtues of its own.

So, let me first derive that 7-term identity.
Applying the second multiplication rule once and then the
first one twice
$$ s(v)s(u)=s(u+v-qvu)=s(u+(v-qvu))
   =s(u)s(v-qvu)=s(u)s(-qvu)s(v)                           $$
we soon come to the 5-term identity\footnote{this `pentagon'
thing leads already its own life, under the banner `Quantum
dilogarithm identity'[FK]}
$$ s(v)s(u)=s(u)s(-qvu)s(v)                                $$
which in turn brings us, again in three
steps\footnote{prior to every step I underline the part
which is going to be treated}, to the 7-term one:
$$ \underline{s(v)}s(u^{-1})\underline{s(u)}s(v)
   =s(u)\underline{s(-qvu)s(v)s(u^{-1})}s(v)\qquad\qquad   $$
$$ \qquad\qquad
   =\underline{s(u)}s(u^{-1})\underline{s(-qvu)s(v)}
   =s(u^{-1})s(v)s(u).                                     $$

One obvious advantage of the 7-term identity, comparing to
the 5-term one and the multiplication rules themselves, is
that we can now produce a closed set
of commutation relations (four nontrivial ones, six in total)
$$ s_2^{+}s_1^{-}s_1^{+}s_2^{+}=s_1^{-}s_2^{+}s_1^{+}\qquad
  \qquad s_2^{-}s_1^{+}s_1^{-}s_2^{-}=s_1^{+}s_2^{-}s_1^{-}$$
$$ s_1^{+}s_2^{+}s_2^{-}s_1^{+}=s_2^{+}s_1^{+}s_2^{-}\qquad
  \qquad s_1^{-}s_2^{-}s_2^{+}s_1^{-}=s_2^{-}s_1^{-}s_2^{+}$$
$$ s_1^{+} s_1^{-} = s_1^{-} s_1^{+}\qquad\qquad
   \qquad s_2^{+} s_2^{-} = s_2^{-} s_2^{+}                $$
involving just four $q$-exponents
$$ s_1^{\pm}= s(u^{\pm 1})\qquad\qquad\qquad
   s_2^{\pm}= s(v^{\pm 1}).                                $$
According to the lattice way of thinking one Weyl pair is
good for a lattice of just two sites. For a longer lattice
one employs a sort of lattice `free field': an algebra where
every `nearest neighbours' pair $w_n,w_{n+1}$ of
its $N$ generators $w_1,w_2,\ldots,w_N$ is like a Weyl pair
$$ w_n w_{n+1}=q^2 w_{n+1}w_n \qquad 1\leq n \leq N-1      $$
while all other pairs just commute
$$ w_m w_n = w_n w_m \qquad |m-n|>1    .                   $$
For $2N$ $q$-exponents available
$$ s_n^{\pm}= s(w_n^{\pm 1})                               $$
$4(N-1)$ nontrivial commutation relations emerge
$$ s_{n+1}^{\pm}s_{n}^{-}s_{n}^{+}s_{n+1}^{\pm}
   =s_{n}^{\mp}s_{n+1}^{\pm}s_{n}^{\pm}\qquad\qquad
   s_{n}^{\mp}s_{n+1}^{-}s_{n+1}^{+}s_{n}^{\mp}
   =s_{n+1}^{\mp}s_{n}^{\mp}s_{n}^{\pm} .                  $$
They are complemented by a bunch of trivial ones
$$ s_m s_n = s_n s_m \qquad |m-n|\neq 1                    $$
where $s_n$ means either $s_n^{+}$ or $s_n^{-}$.
Meet a brand new discrete group. Indeed, we can
now dispose of the free field and regard $s$'s as
just generators obeying only the above set of
commutation relations.

Of course, the crucial question is whether or not the 7-term
identity is all we really want to know about
the $q$-exponent. 
Apparently it is, at least as far as the Lattice CFT [FV] is
concerned. First come the braids. The elements
$$ b_n=s_n^{+}s_n^{-}                                      $$
$$ b_m b_n = b_n b_m \qquad |m-n|>1                        $$
prove to obey the Artin's commutation relations:
$$ b_{n} b_{n+1} b_{n} 
   =s_{n}^{-}
   \underline{s_{n}^{+}s_{n+1}^{+}s_{n+1}^{-}s_{n}^{+}}
   s_{n}^{-} 
   =\underline{s_{n}^{-}s_{n+1}^{+}s_{n}^{+}}
   s_{n+1}^{-}s_{n}^{-}\qquad\qquad\qquad                  $$
$$ =s_{n+1}^{+}s_{n}^{+}
   \underline{s_{n}^{-}s_{n+1}^{+}s_{n+1}^{-}s_{n}^{-}}    $$
$$ \qquad\qquad\qquad=s_{n+1}^{+}
   \underline{s_{n}^{+}s_{n+1}^{-}s_{n}^{-}}s_{n+1}^{+}
   =s_{n+1}^{+}s_{n+1}^{-}s_{n}^{+}
   s_{n}^{-}s_{n+1}^{-}s_{n+1}^{+}
   = b_{n+1} b_{n} b_{n+1}     .                           $$
This is indeed the braid group $B_{N+1}$.
It is however remains to see what the `twisted' set-up
$$ \varsigma_n=s_n^{-}s_{n+1}^{+}                          $$
can do. Fortunately, it delivers: 
$$ \varsigma_{n+1}\varsigma_{n-1}\varsigma_{n}\varsigma_{n+1}
   =s_{n+1}^{-}\underline{s_{n+2}^{+}}s_{n-1}^{-}s_{n}^{+}
   s_{n}^{-}\underline{s_{n+1}^{+}s_{n+1}^{-}s_{n+2}^{+}}
   \qquad\qquad\qquad$$
$$ =\underline{s_{n+1}^{-}}s_{n-1}^{-}
   \underline{s_{n}^{+}s_{n}^{-}s_{n+1}^{-}}
   s_{n+2}^{+}s_{n+1}^{+}                                  $$
$$ \qquad\qquad\qquad
   =s_{n-1}^{-}s_{n}^{+}s_{n+1}^{-}s_{n}^{-}     
   s_{n+2}^{+}s_{n+1}^{+}
   =\varsigma_{n-1}\varsigma_{n+1}\varsigma_{n}.           $$
Similarly,
$$ \varsigma_{n-1}\varsigma_{n}\varsigma_{n+1}\varsigma_{n-1}
   =\underline{s_{n-1}^{-}s_{n}^{+}s_{n}^{-}}s_{n+1}^{+}
   s_{n+1}^{-}s_{n+2}^{+}\underline{s_{n-1}^{-}}s_{n}^{+}
   \qquad\qquad\qquad                                      $$
$$ =s_{n}^{-}s_{n-1}^{-}
   \underline{s_{n}^{+}s_{n+1}^{+}s_{n+1}^{-}}
   s_{n+2}^{+}\underline{s_{n}^{+}}                        $$
$$ \qquad\qquad\qquad
   =s_{n}^{-}s_{n-1}^{-}s_{n+1}^{+}s_{n}^{+}
   s_{n+1}^{-}s_{n+2}^{+}
   =\varsigma_{n}\varsigma_{n-1}\varsigma_{n+1} .          $$
So, we end up with yet another group and there is a good
reason to call (a cyclic version of) this one a `lattice
Virasoro algebra' or maybe a `discrete conformal
group'. This issue, as well as that of
YangBaxterization, will be discussed in detail elsewhere.

Anyway, the group
$$ \varsigma_{n+1}\varsigma_{n-1}\varsigma_{n}\varsigma_{n+1}
   =\varsigma_{n-1}\varsigma_{n+1}\varsigma_{n}            $$
$$ \varsigma_{n-1}\varsigma_{n}\varsigma_{n+1}\varsigma_{n-1}
   =\varsigma_{n}\varsigma_{n-1}\varsigma_{n+1}            $$

$$ \varsigma_m\varsigma_n=\varsigma_n\varsigma_m
   \qquad |m-n|>2                                          $$
seems to be
the most valuable outcome of those $q$-manipulations. It
looks like a close relative to the braid group
$$ b_{n} b_{n+1} b_{n} = b_{n+1} b_{n} b_{n+1}             $$
$$ b_m b_n = b_n b_m \qquad |m-n|>1                        $$
for despite of their obvious differences they still share
some key features. One striking similarity between them
is how a single generator goes through long
enough `ordered' words:
$$ (b_m b_{m+1} \ldots b_n)b_k
  =b_m\ldots (b_k b_{k+1} b_k) \ldots b_n
   \qquad\qquad\qquad                                      $$
$$ \qquad\qquad\qquad
   =b_m\ldots (b_{k+1} b_k b_{k+1})\ldots b_n
   =b_{k+1}(b_m b_{m+1} \ldots b_n)                        $$
$$ (\varsigma_{m}\varsigma_{m+1}\ldots \varsigma_{n})
   \varsigma_{k}=\varsigma_{m}\ldots\varsigma_{k-1}
   (\varsigma_{k}\varsigma_{k+1}\varsigma_{k+2}
   \varsigma_{k})\ldots \varsigma_{n}
   \qquad\qquad\qquad\qquad                                $$
$$ =\varsigma_{m}\ldots\varsigma_{k-1}
   (\varsigma_{k+1}\varsigma_{k}\varsigma_{k+2})
   \ldots \varsigma_{n}=\varsigma_{m}\ldots(\varsigma_{k-1}
   \varsigma_{k+1}\varsigma_{k})\varsigma_{k+2}
   \ldots \varsigma_{n}                                    $$
$$ \qquad\qquad\qquad\qquad=\varsigma_{m}
   \ldots(\varsigma_{k+1}\varsigma_{k-1}
   \varsigma_{k}\varsigma_{k+1})\varsigma_{k+2}
   \ldots \varsigma_{n}=\varsigma_{k+1}
   (\varsigma_{m}\varsigma_{m+1}\ldots \varsigma_{n}).     $$
Of course, the similarity can not remain this literal for
reversely ordered words but it appears no less amusing.
While in the braid group this is again a one-step translation
$$ (b_n b_{n-1} \ldots b_m)b_{k+1}
   =b_k(b_n b_{n-1} \ldots b_m)  ,                         $$
in $\varsigma$'s it is a translation by two steps at once: 
$$ (\varsigma_{n}\varsigma_{n-1}\ldots \varsigma_{m})
   \varsigma_{k+1}=\varsigma_{n}\ldots\varsigma_{k+1}
   (\varsigma_{k}\varsigma_{k-1}\varsigma_{k+1})\ldots
   \varsigma_{m}\qquad\qquad\qquad\qquad                   $$
$$ =\varsigma_{n}\ldots\varsigma_{k+1}(\varsigma_{k-1}
   \varsigma_{k}\varsigma_{k+1}\varsigma_{k-1})\ldots
   \varsigma_{m}
   =\varsigma_{n}\ldots(\varsigma_{k+1}\varsigma_{k-1}
   \varsigma_{k}\varsigma_{k+1})\varsigma_{k-1}\ldots
   \varsigma_{m}                                           $$
$$ \qquad\qquad\qquad\qquad
   =\varsigma_{n}\ldots(\varsigma_{k-1}\varsigma_{k+1}
   \varsigma_{k})\varsigma_{k-1}\ldots\varsigma_{m}
   =\varsigma_{k-1}
   (\varsigma_{n}\varsigma_{n-1}\ldots \varsigma_{m}).     $$
Not really proving anything, these simple tests at least
give a hope that the new group is only about as `large' as
the braid group. And if this is indeed true, it might find
a spectrum of applications reaching far beyond our
modest Lattice CFT.
\bigskip
\begin{quote}
{\bf Acknowledgements.}
I would like to thank A.~Izergin, L.~Faddeev, V.~Fateev,
R.~Kashaev, J.-M.~Maillet, A.~Neveu and A.~Reyman
for stimulating discussions.
\end{quote}
\section*{\normalsize\bf References}

\begin{itemize}
  \item[{[GR]}] G. Gasper and M. Rahman,
Encyclopedia of Mathematics and its Applications 35,\\
Cambridge University Press, 1990.
  \item[{[FV]}] L. Faddeev and A. Yu. Volkov,
Phys. Lett. B315 (1993) 311.
  \item[{[FK]}] L. Faddeev and R. M. Kashaev,
Modern Phys. Lett. A9 (1994) 427.
\end{itemize}
\end{document}